\documentclass[preprint,11pt]{elsarticle}
\usepackage{amssymb,amsmath}
\usepackage{graphics}
\biboptions{numbers,sort&compress}
\newtheorem{theorem}{Theorem}[section]
\newtheorem{lemma}{Lemma}[section]
\newtheorem{remark}{Remark}[section]
\newtheorem{cor}{Corollary}
\begin{document}

\begin{frontmatter}

\title{ Remarks on Bounds of Normalized Laplacian Eigenvalues of  Graphs}


\author{ Emina I. Milovanovi\' c, Igor \v Z. Milovanovi\' c}

\address{Faculty of Electronic Engineering, A. Medvedeva 14, P.O.Box 73, 18000 Ni\v s, Serbia}

\begin{abstract}
Let $G$ be a connected undirected graph with $n$, $n\ge 3$,  vertices and $m$ edges. Denote by $\rho_1 \ge \rho_2 \ge \cdots > \rho_n =0$ the normalized  Laplacian eigenvalues of  $G$. Upper and lower bounds of $\rho_i$, $i=1,2,\ldots , n-1$, are determined in terms of $n$ and general Randi\' c index, $R_{-1}$.
\end{abstract}


\end{frontmatter}

\section{Introduction and preliminaries}

Let $G=(V,E)$, $V=\{1,2, \ldots , n\}$, be $n$-vertex   undirected connected graph, with  $m$ edges and vertex degree sequence  $d_1\ge d_2 \ge \cdots \ge d_n>0$, $d_i=d(i)$. Denote by $\mathbf A$ adjacency matrix of $G$ and by $\mathbf D =\mbox{diag}\{d_1,d_2,\ldots , d_m\}$ a diagonal matrix of its vertex degrees.Then, $\mathbf L =\mathbf D -\mathbf A$ is Laplacian matrix of $G$. Since it is assumed that $G$ is connected, matrix $\mathbf D$ is nonsingular, and thus, $\mathbf D^{-1}$ always exists. Matrix $\mathbf L^* =\mathbf D^{-1/2}\mathbf L \mathbf D^{-1/2} = \mathbf I -\mathbf D^{-1/2}\mathbf  A \mathbf D^{-1/2}$ is called normalized Laplacian matrix of  $G$. Eigenvalues of $\mathbf L^*$, $\rho_1\ge \rho_2 \ge \cdots > \rho_n = 0$,  are normalized Laplacian eigenvalues of graph $G$ (see \cite{c3}). Well known properties of these eigenvalues are \cite{c11}
\begin{equation}
\sum_{i=1}^{n-1} \rho_i =n \qquad\mbox{ and } \quad \sum_{i=1}^{n-1} \rho_i^2 =n+2R_{-1},
\label{l1.1}
\end{equation}
where $R_{-1}$ is general Randi\' c index \cite{c1,c2,c8,c11}. Well known  inequalities valid for  $\rho_1$ and $\rho_{n-1}$ are \cite{c3}
\begin{equation}
\rho_1 \ge \frac{n}{n-1} \qquad\mbox{ and }\quad \rho_{n-1} \le \frac{n}{n-1},
\label{l1.2}
\end{equation}
with equality holding if and only if $G\cong K_n$. In \cite{c1,c2,c3,c4,c5,c6,c9,c10} several upper/lower bounds for $\rho_1$ and $\rho_{n-1}$ were reported. In the present paper we consider upper and lower bounds for $\rho_i$, $i=1,2,\ldots ,n-1$.

In the text that follows we first recall some results from spectral graph theory and polynomials needed for our work.

Shi \cite{c9} (see also \cite{c4}) proved the following result for general Randi\' c index.
\begin{lemma}
\label{le1.1}{\rm \cite{c9}}
Let $G$ be a connected undirected  graph with $n$ vertices and $m$ edges. Then
\begin{equation}
\frac{n}{2d_1} \le R_{-1} \le \frac{n}{2d_n} .
\label{l1.3}
\end{equation}
Equality holds if and only if $G$ is a regular graph.
\end{lemma}

Lupas \cite{c7} considered a class of real polynomials ${\cal P}_n(a_1,a_2)$, i.e. the polynomials of type $P_n(x) =x_n+a_1x^{n-1}+a_2x^{n-2}+b_3x^{n-3} +\cdots +b_n$, where $a_1$ and $a_2$ are fixed real numbers. Denote by $x_1\ge x_2\ge \cdots \ge x_n$ zeros of polynomial $P_n(x)$, whereby
\begin{equation}
\bar x = \frac 1n \sum_{i=1}^n x_i \qquad\mbox{ and }\quad \Delta = n \sum_{i=1}^n x_i^2 -\left(\sum_{i=1}^n x_i\right)^2
\label{l1.4}
\end{equation}
The following bounds for zeros of polynomial $P_n(x)$ were proved in \cite{c7}.

\begin{lemma}
\label{le1.2}
{\rm \cite{c7}} If $P_n(x)\in {\cal P}_n(a_1,a_2)$ then
\begin{eqnarray}
&&\bar x +\frac 1n \sqrt{\frac{\Delta}{n-1}} \le x_1 \le \bar x +\frac 1n \sqrt{(n-1)\Delta} , \nonumber\\
&& \bar x -\frac 1n \sqrt{\frac{i-1}{n-i+1}\Delta} \le x_i \le \bar x +\frac 1n \sqrt{\frac{n-i}{i}\Delta}, \qquad 2 \le i \le n-1 \label{l1.5}\\
&& \bar x -\frac 1n \sqrt{(n-1)\Delta}  \le x_n \le \bar x - \frac 1n \sqrt{\frac{\Delta}{n-1}} . \nonumber
\end{eqnarray}
\end{lemma}

\section{Main results}

In the following theorem we determine upper and lower bounds  for normalized Laplacian eigenvalues, $\rho_i$, $i=1,2,\ldots , n-1$ in terms of $n$ and $R_{-1}$.

\begin{theorem}
Let $G$ be a connected undirected  graph with $n$, $n\ge 3$, vertices and $m$ edges. Then

\begin{align}
\label{l2.1}
\begin{split}
& \frac{n}{n-1} + \frac{1}{n-1}\sqrt{\frac{2(n-1)R_{-1}-n}{n-2}} \le \rho_1 \le \\
&\le \frac{n}{n-1} +\frac{1}{n-1}\sqrt{(n-2)(2(n-1)R_{-1}-n)},
\end{split}
\end{align}
\begin{eqnarray*}
&& \frac{n}{n-1} -\frac{1}{n-1} \sqrt{\frac{i-1}{n-i}(2(n-1)R_{-1} -n)} \le \rho_i \le\\
&\le& \frac{n}{n-1} +\frac{1}{n-1} \sqrt{\frac{n-i-1}{i}(2(n-1)R_{-1} -n)}, \, 2\le i\le n-2
\end{eqnarray*}

\begin{align}
\begin{split}
\frac{n}{n-1} - \frac{1}{n-1}&\sqrt{(n-2)(2(n-1)R_{-1} -n)} \le \rho_{n-1} \le \\ &\le \frac{n}{n-1} - \frac{1}{n-1} \sqrt{\frac{2(n-1)R_{-1}-n}{n-2}} .
\end{split}
\label{l2.2}
\end{align}
In all three cases equalities hold if and only if $G\cong K_n$.
\end{theorem}

\noindent{\bf Proof.} The characteristic polynomial of the normalized Laplacian matrix $\mathbf L^*$  of graph $G$ is 
$$
\bar \varphi_{n}(x) = x\varphi_{n-1}(x) = x\left(x^{n-1} +a_1x^{n-2} +a_2x^{n-3} +b_3x^{n-4} +\cdots + b_{n-1}\right),
$$
whereby
$$
a_1 = -\sum_{i=1}^{n-1} \rho_i = -n \quad\mbox{ and } a_2 =\frac 12\left( \left( \sum_{i=1}^{n-1} \rho_i\right)^2 -\sum_{i=1}^{n-1}\rho_i^2\right) = \frac{n(n-1)}{2} - R_{-1}.
$$
This means that polynomial $\varphi_{n-1}(x)$ belongs to a class of polynomial \\ ${\cal P}_{n-1}(-n,\frac{n(n-1)}{2} -R_{-1})$. According to (\ref{l1.4}) for the zeros, $\rho_1\ge \rho_2 \ge \cdots \ge \rho_{n-1} >0$,  of polynomial $\varphi_{n-1}(x)$ the following equalities are valid
\begin{align}
\label{l2.3}
\begin{split}
\bar x &=\frac{1}{n-1}\sum_{i=1}^{n-1} \rho_i =\frac{n}{n-1},  \\ \Delta &= (n-1)\sum_{i=1}^{n-1}\rho_i^2  -\left(\sum_{i=1}^{n-1}\rho_i\right)^2 = 2(n-1)R_{-1} -n .
 \end{split}
\end{align}
Now according to (\ref{l2.3}) and from (\ref{l1.5}) for $n:=n-1$, $x_i:=\rho_i$, $i=1,2, \ldots , n-1$ we obtain the required result.
\begin{flushright}
    $\blacksquare$
\end{flushright}

\begin{remark}
    \label{r2.1}
    Left-side inequality in (\ref{l2.1}) and right-side inequality in (\ref{l2.2}) were proved in \cite{c10}.
    \end{remark}

\begin{cor}
\label{cor2.1}
Let $G$ be a connected undirected  graph with $n$, $n\ge 3$, vertices and $m$ edges. Then
\begin{align}
\label{l2.4}
\begin{split}
\frac{n}{n-1} +\frac{1}{n-1} &\sqrt{\frac{n(n-1-d_1}{(n-2)d_1}} \le \rho_1 \le\\
& \le \frac{n}{n-1} +\frac{1}{n-1}\sqrt{\frac{n(n-2)(n-1-d_n)}{d_n}}
 \end{split}
\end{align}
\begin{eqnarray*}
&& \frac{n}{n-1} - \frac{1}{n-1}\sqrt{\frac{n(i-1)(n-1-d_n)}{(n-i)d_n}} \le \rho_i \le\\ &&\le \frac{n}{n-1}  +\frac{1}{n-1} \sqrt{\frac{n(n-i-1)(n-1-d_n)}{id_n}}, \, 2\le i\le n-2
\end{eqnarray*}
\begin{align}
\label{l2.5}
\begin{split}
\frac{n}{n-1} - \frac{1}{n-1} &\sqrt{\frac{n(n-2)(n-1-d_n)}{d_n}} \le \rho_{n-1} \le\\
& \le \frac{n}{n-1} - \frac{1}{n-1}\sqrt{\frac{n(n-1-d_1)}{(n-2)d_1}}.
 \end{split}
\end{align}
In all three cases equalities hold if and only if $G\cong K_n$.
\end{cor}

\begin{remark}
\label{r2.2}
Let us note that left-side inequalities in (\ref{l2.1}) and (\ref{l2.4}) as well as right-side inequalities in (\ref{l2.2}) and (\ref{l2.5}) are stronger than those in (\ref{l1.2}).
\end{remark}

\begin{remark}
\label{r2.3}
Based on the left inequality in (\ref{l2.1}) and the right inequality in (\ref{l2.2}) the following inequalities are obtained, respectively
$$
R_{-1} \le \frac 12 (n-1)(n-2) \left( \rho_1 - \frac{n}{n-1}\right)^2 +\frac{n}{2(n-1)}
$$
and
$$
R_{-1} \le \frac 12 (n-1)(n-2)\left( \frac{n}{n-1} -\rho_{n-1}\right)^2 +\frac{n}{2(n-1)},
$$
which were proved in \rm \cite{c10} [Theorem 1.1].
\end{remark}
\begin{cor}
\label{cor2.2}
Let $G$ be a connected undirected  graph with $n$, $n\ge 3$, vertices and $m$ edges. Then
\begin{equation}
R_{-1} \ge \frac{n-1}{2(n-2)}\left( \rho_1 -\frac{n}{n-1}\right)^2 +\frac{n}{2(n-1)}
\label{l2.6}
\end{equation}
and
\begin{equation}
R_{-1} \ge \frac{n-1}{2(n-2)}\left(\frac{n}{n-1} -\rho_{n-1}\right)^2 +\frac{n}{2(n-1)}
\label{l2.7}
\end{equation}
Equalities hold if and only if $G\cong K_n$.
\end{cor}

\begin{remark}
\label{r2.4}
The inequality (\ref{l2.6}) was proved in \cite{c10} [Theorem 1.2], under slightly stronger condition. Also, it is not difficult to see that inequalities (\ref{l2.6}) and (\ref{l2.7}) are stronger than inequality $R_{-1} \ge \frac{n}{2(n-1)}$ proved in \cite{c6}.
\end{remark}

\end{document}